\documentclass[11pt,reqno]{amsart}
\usepackage[utf8]{inputenc}
\usepackage[T1]{fontenc}
\usepackage{lmodern}
\usepackage{amssymb}
\usepackage{tikz}
\usepackage[margin=1in]{geometry}
\usepackage[colorlinks,linkcolor=blue,citecolor=red,urlcolor=blue]{hyperref}

\theoremstyle{plain}
\newtheorem{theorem}{Theorem}[section]
\newtheorem{lemma}[theorem]{Lemma}
\newtheorem{proposition}[theorem]{Proposition}
\theoremstyle{remark}
\newtheorem{remark}[theorem]{Remark}
\numberwithin{equation}{section}

\newcommand{\PP}{\mathbb P}
\newcommand{\EE}{\mathbb E}
\newcommand{\ZZ}{\mathbb Z}
\newcommand{\NN}{\mathbb N}
\newcommand{\RR}{\mathbb R}
\newcommand{\FF}{\mathcal F}
\newcommand{\ind}{\mathbf 1}

\author[Shuo Qin]{Shuo Qin}
\address{Beijing Institute of Mathematical Sciences and Applications, and Yau Mathematical Sciences Center, Tsinghua University}
\email{qinshuo@bimsa.cn}

\title[Transience of overlapping balanced excited walks]{Transience of overlapping balanced excited random walks in dimensions at least three}
\date{}

\begin{document}
\begin{abstract}
We study a walk on $\ZZ^3$ whose first departure from each vertex is a
simple random walk step in the $xy$-plane and whose later departures are
simple random walk steps in the $yz$-plane. We prove that it is at
distance $n^{1/2+o(1)}$ from the origin at time $n$, almost surely. The same escape rate holds for two further walks in dimensions three and four. Together with earlier results and applications of a trace criterion, this proves transience for all overlapping balanced excited random walks in dimensions at least three.
\end{abstract}
\maketitle

\section{Introduction}

A simple way for a random walk to depend on its past is to distinguish first and subsequent departures from a vertex. In the excited random
walk introduced by Benjamini and Wilson \cite{MR1987097}, the first
departure has a bias in a fixed direction, while later departures
follow the simple random walk rule. More generally, the transition
probabilities may depend on the number of previous visits to the
current vertex. We refer to Kosygina and Zerner \cite{MR3097419} for a
survey. 

More recently, Kosygina and Peterson \cite{MR3689974} studied one-dimensional excited
walks with Markovian cookie stacks and obtained criteria for recurrence
and transience. Recent work has also considered models in which only the first
departure from each vertex uses a different transition rule.
Le and Nguyen \cite{LN2026once} established recurrence and transience
criteria for once-excited walks on trees of polynomial growth.
Bou-Rabee, Nguyen, Peres and Singh \cite{BRNPS2026} proved a shape
theorem for centrally excited walks on $\ZZ^d$, $d\ge2$, whose first departure from each nonzero vertex has a drift toward the origin:
the visited set is asymptotically a Euclidean ball.

We study balanced excited walks, for which both transition laws have mean zero but may allow motion in different coordinate directions. 

\subsection{The model and main results}
Let $e_1,\ldots,e_d$ denote the standard coordinate vectors of $\RR^d$.
For $1\le d_1,d_2\le d$ with $d_1+d_2\ge d$, the balanced excited walk
$M(d,d_1,d_2)$ starts at the origin of $\ZZ^d$. On its first departure
from a vertex its increment is uniform on $\{\pm e_i:1\le i\le d_1\}$.
On every later departure from that vertex the increment is uniform on
$\{\pm e_i:d-d_2+1\le i\le d\}$. Both laws have mean zero. Their coordinate
sets together span $\RR^d$ and overlap if $d_1+d_2>d$.
We call the walk transient if almost surely every vertex is visited only
finitely often.

Benjamini, Kozma and Schapira \cite{MR2788390} introduced the
non-overlapping case $d_1+d_2=d$ and proved transience for $M(4,2,2)$.
In the other non-overlapping cases with $d\ge4$, projection onto the coordinate set with dimension at least three implies transience. Peres, Schapira and Sousi \cite{MR3531697} proved transience
for $M(3,1,2)$ by a martingale defocusing argument. Camarena, Panizo and
Ram\'irez \cite{Camarena2020} introduced the overlapping case and proved
transience in dimensions at least four except for $M(4,3,2)$. In two dimensions, Angel, Holmes and Ram\'irez \cite{MR4597323} obtained
upper and lower bounds for the range of $M(2,1,1)$, whose recurrence
remains open. Qin \cite{Qin2026M212} proved recurrence and obtained the
range asymptotics of $M(2,1,2)$. Peres and
Qin \cite{PeresQin2026Columns} proved transience for $M(3,2,1)$ and
$M(3,3,1)$ by estimating occupation times in vertical columns. 

We prove transience for the remaining cases in dimensions $d\geq 3$.

\begin{theorem}\label{thm-m322}
Let $S$ be the $M(3,2,2)$ walk. There is a constant $c>0$ such that for
every fixed $v\in\ZZ^3$ there is $C_v<\infty$ with
\begin{equation}\label{eq-tail}
 \sum_{k\ge n}\PP(S_k=v)\le C_v\exp\{-c(\log n)^{1/4}\},\qquad n\ge 3 .
\end{equation}
Moreover, with $|\cdot|$ denoting the Euclidean norm,
\begin{equation}\label{eq-escape}
 \lim_{n\to\infty}\frac{\log(1+|S_n-v|)}{\log n}=\frac12
 \qquad\text{almost surely}.
\end{equation}
\end{theorem}

\begin{theorem}\label{thm-m3-432}
Let $S$ be either the $M(3,3,2)$ walk or the $M(4,3,2)$ walk. Then
the occupation bound \eqref{eq-tail} and the escape limit \eqref{eq-escape}
hold for every fixed vertex $v$ of the corresponding lattice.
\end{theorem}

 The walk $M(3,3,3)$ is ordinary three-dimensional simple random walk. For $d\geq 3$, the remaining overlapping cases are $M(3,1,3)$ and $M(3,2,3)$.
\begin{proposition}\label{prop-all3models}
The walks $M(3,1,3)$ and $M(3,2,3)$ are transient. 
\end{proposition}
Consequently, $M(d,d_1,d_2)$ is transient for every $d\ge3$.

\subsection{Proof strategy}
We first outline the proof of Theorem \ref{thm-m322} for $M(3,2,2)$.
By translation we may take the target vertex to be the origin. Write
$\PP_{z_0}$ for the law of the walk started at $z_0$ with no previous
visits, and define
\[
 \tau_s=\inf\{n\ge0:|S_n|\ge s\},\qquad
 d_R=R\exp\{-(\log R)^{3/4}\}.
\]
The main step is Proposition \ref{prop-cross}, which bounds the
probability of returning close to the origin after reaching distance
$R$. For all sufficiently large $R$, uniformly in $z_0$, it gives
\[
 \PP_{z_0}\bigl(\exists\,n\in[\tau_R,\tau_{2R}): |S_n|\le d_R\bigr)
 \le C\exp\{-c(\log R)^{1/4}\}.
\]
This bound is summable over dyadic $R$. Since $\log d_R\sim\log R$,
once we know that $\tau_R=R^{2+o(1)}$ almost surely along these radii,
Borel--Cantelli gives the escape rate in \eqref{eq-escape}. To obtain
\eqref{eq-tail}, we also bound the expected number of visits following
a first return to the origin, sum over dyadic radii, and use exit-time
estimates to pass from $\tau_R$ to a deterministic time.

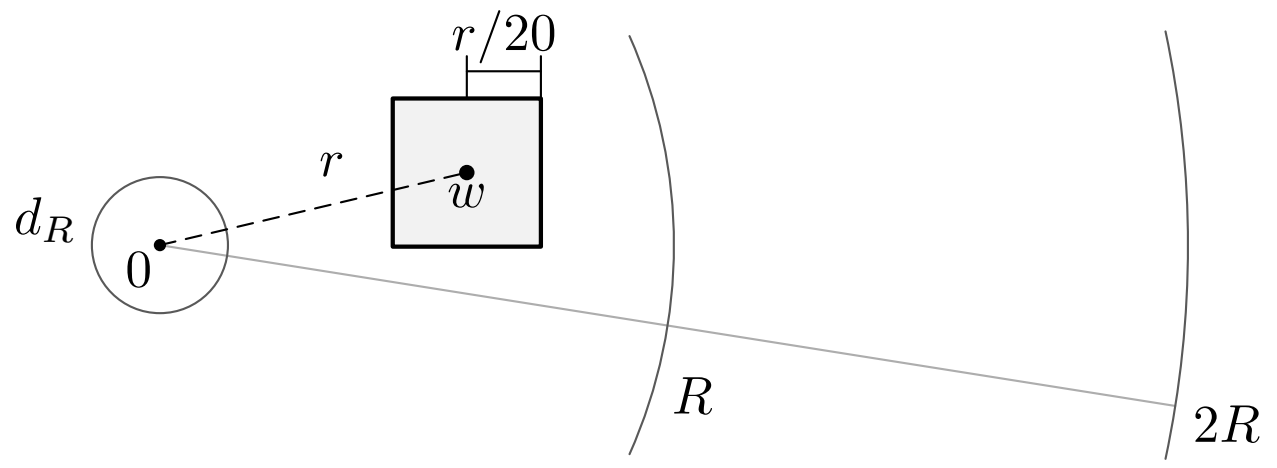
\begin{figure}[!ht]
\centering
\begin{tikzpicture}[font=\small,line cap=round,line join=round]
  \draw[black!32,thin] (0,0) -- (-9:6.8);
  \draw[black!65] (0,0) circle[radius=.45];
  \draw[black!65] (-24:3.4) arc[start angle=-24,end angle=24,radius=3.4];
  \draw[black!65] (-12:6.8) arc[start angle=-12,end angle=12,radius=6.8];
  \node[left] at (158:.45) {$d_R$};
  \node[right] at (-17:3.4) {$R$};
  \node[right] at (-10:6.8) {$2R$};
  \coordinate (w) at (2.03,.48);
  \filldraw[fill=black!5,draw=black,thick] (1.54,-.01) rectangle (2.52,.97);
  \draw[densely dashed] (0,0) -- (w) node[pos=.56,above=1.5pt] {$r$};
  \fill (0,0) circle[radius=1.15pt] node[below left,inner sep=1.3pt] {$0$};
  \fill (w) circle[radius=1.45pt] node[below,inner sep=2pt] {$w$};
  \draw[thin] (2.03,.97) -- (2.03,1.25);
  \draw[thin] (2.52,.97) -- (2.52,1.25);
  \draw (2.03,1.15) -- (2.52,1.15)
    node[midway,above,inner sep=1.3pt] {$r/20$};
\end{tikzpicture}
\caption{The spheres of radii $d_R,R,2R$ about the origin and a cube centered at the walk's position $w$ with half-side $r/20$ where $r = |w|$.}
\label{fig-cube}
\end{figure}

Fix an outer scale $R$. We construct a positive function $V = V_R$,
comparable to $|v|^{-\kappa}$ with $\kappa \asymp (\log R)^{-1/2}$ and modified by an angular factor. We sample the walk at successive exits from cubes centered at its current position. Each cube has half-side equal to a small fixed multiple of the distance of its center from the origin, see
Figure \ref{fig-cube} for an illustration. Using the angular correction, we show that the conditional expectation of $V$ does not increase over a cube exit whenever the projected past is sufficiently sparse.

More specifically, under the law of $M(3,2,2)$, the projection $(x,y,z)\mapsto(x+z,y)$ is a planar simple random walk: its next increment is uniform on the four planar directions, conditional on the history. We group departures by their starting sites and use planar Green-function estimates to control the expected visit counts. At sites whose projections were unvisited before the cube was chosen, bounds on the relative numbers of first and later departures give a negative contribution to the expected change of $V$. The possible positive contribution from the remaining sites is controlled by their total expected occupation, which is small when the projected past is sparse.

We then iterate the cube-exit inequality, stopping when the walk
approaches the inner ball, leaves the outer ball, or the sparsity
condition fails. The resulting supermartingale bounds the probability of approaching the origin by $C(d_R/R)^\kappa$, while a planar range estimate bounds the probability that sparsity fails at any sampling time; the choice $d_R = R\exp\{-(\log R)^{3/4}\}$ balances these two bounds at $\exp\{-c(\log R)^{1/4}\}$. The cubes are small enough that an inward crossing requires one of these two events. Combining the
estimates proves Proposition \ref{prop-cross}.

Section \ref{sec-transfer} adapts the cube-exit argument and the same
test function to $M(3,3,2)$ and $M(4,3,2)$, after projecting the
four-dimensional walk onto $\ZZ^3$. Finally,
Section \ref{sectrace} proves transience of $M(3,1,3)$ and
$M(3,2,3)$ by applying the trace condition of Peres, Popov and
Sousi \cite{MR3167134} to three-step increments.

Throughout the paper, $C$ and $c$ denote positive constants that may change from line to line, and $f \asymp g$ means $cg \le f \le Cg$.

\section{Green-function and sparse-range estimates}

We need two estimates for the planar projection used in
Section \ref{sec-m322}. The first bounds the expected time spent in a
finite set before leaving a disk. The second shows that the projected
past occupies only a small fraction of each disk needed in the cube
argument.

Write $B(x,r):=\{u\in\ZZ^2:|u-x|<r\}$ and $B(r):=B(0,r)$.
Throughout this section, $U$ is planar simple random walk with law
$\mathrm P_x$ and expectation $\mathrm E_x$ when started at $x$.
For a finite set $A$, let $T_{A^c}$ be its first exit time and $T_u$
the first visit to $u$, including time zero. The killed Green function
\[
 G_A(x,u):=\mathrm E_x\sum_{n<T_{A^c}}\ind_{\{U_n=u\}}
\]
is the expected number of visits to $u$ before exit. In particular, $G_A$ is increasing in $A$ with respect to the set inclusion. By \cite[Lemma 4.6.1]{MR2677157}, it is symmetric and related to the hitting probability by
\begin{equation}\label{eq-hitquot}
 \mathrm P_x(T_u<T_{A^c})=\frac{G_A(x,u)}{G_A(u,u)},\qquad x,u\in A.
\end{equation}

\begin{lemma}\label{lem-green}
Let $A=B(w,r)$ with $r\ge2$. Then one has the following diagonal bounds:
\[
 G_A(u,u)\le C\log r\quad (u\in A),\qquad
 G_A(u,u)\ge c\log r\quad (|u-w|\le r/2).
\]
For every $x\in A$ and $E\subseteq A$ of size $m$,
\begin{equation}\label{eq-setgreen}
 \sum_{u\in E}G_A(x,u)
 \le C m\big[1+\log_+(r/\sqrt{1+m})\big],\qquad
 \log_+s:=\max(0,\log s).
\end{equation}
\end{lemma}
\begin{proof}
For $u\in A$ we have $A\subseteq B(u,2r)$, and
$B(u,r/2)\subseteq A$ when $|u-w|\le r/2$. Also note that $G_{B(s)}(0,0)\asymp\log s$ for $s\ge2$, by \cite[Proposition 6.3.2]{MR2677157} (for simple random walk, the norm
$\mathcal J$ used there is the Euclidean norm). These inclusions imply the desired diagonal upper bound for all $r\geq 2$, and lower bound for $r\geq 4$.  The lower bound for
$2\le r<4$ follows from $G_A(u,u)\ge1$.

For $x\ne u$, symmetry and
\cite[Proposition 6.3.5]{MR2677157} give
\[
 G_A(x,u)\le G_{B(u,2r)}(u,x)
       \le C\big[1+\log_+(r/(1+|x-u|))\big].
\]
The diagonal upper bound above gives the same estimate at $x=u$.
To sum this bound over $E$, we use Tonelli's theorem and the fact that
$|E\cap B(x,s)|\le C\min(m,s^2)$ for $s\ge1$. Thus
\[
\begin{aligned}
 \sum_{u\in E}\log_+\frac{r}{1+|x-u|}
 &\le  \sum_{u\in E} \int_1^r \mathbf{1}_{\{|x-u|<s\}} \, \frac{ds}{s} \le C\int_1^r\min(m,s^2)\,\frac{ds}{s} \\
 & \le C m\big[1+\log_+(r/\sqrt{1+m})\big].
\end{aligned}
\]
The last inequality follows by splitting the integral at
$\sqrt{1+m}$ when this lies below $r$. Adding the constant terms gives
\eqref{eq-setgreen}.
\end{proof}

The next estimate concerns the number of distinct visited sites, rather
than the time spent at them. Peres, Schapira and Sousi
\cite[Lemma 3.5]{MR3531697} estimated the number of sites in a disk
visited by a planar walk up to a deterministic time. Here we need a
version up to an exit time, uniform over the radii that occur in the
cube argument. For the rest of the paper, write $L=\log R$ for simplicity, so that $ d_R=R e^{-L^{3/4}}$.
\begin{lemma}\label{lem-sparse}
Let $\mathcal U_R$ be the range of $U$ before its first exit from the ball 
$B(32R)$. For every fixed $\delta>0$ and all sufficiently large $R$,
\begin{equation}\label{eq-sparse}
 \sup_{x\in B(32R)}\mathrm P_x\!\left(
  \exists r\in[2d_R,2R]:|\mathcal U_R\cap B(3r)|>\delta r^2\right)
 \le C_\delta e^{-c_\delta L^{1/4}} .
\end{equation}
\end{lemma}
\begin{proof}
We first fix a radius $s\in[2d_R,12R]$ and estimate
$J_s:=|\mathcal U_R\cap B(s)|$. By Lemma \ref{lem-green} with $A=B(32R)$, every $u\in B(s)$ has $G_A(u,u)\ge cL$. Since $\mathrm{E}_x J_s = \sum_{u \in B(s)} \mathrm{P}_x(T_u < T_{A^c})$, summing
\eqref{eq-hitquot} over $u\in B(s)$ and using \eqref{eq-setgreen}
with $|B(s)|\asymp s^2$ gives, for all large $R$,
\[
 \sup_{x\in A}\mathrm E_x J_s
 \le C s^2\,\frac{1+\log(32R/s)}{L}
 \le \tilde{C} s^2 L^{-1/4},
\]
where $\tilde{C}$ is a positive constant, and we used $s\ge2d_R$ in the last inequality.

We turn this first-moment bound into an exponential tail by restarting
the walk at successive discovery times. Put
$k=\lceil 2\tilde{C}s^2 L^{-1/4}\rceil$. Markov's inequality shows that, from
any starting point in $A$, the probability of visiting $k$ distinct
sites of $B(s)$ before exit is at most $1/2$. After the $jk$-th
discovery, the number of additional sites discovered is bounded by the
number of distinct sites visited by the restarted walk. The strong
Markov property therefore gives the same upper bound $1/2$ for
making a further $k$ discoveries. Iteration yields
\[
 \mathrm P_x(J_s\ge jk)\le2^{-j},\qquad j\ge1.
\]
For each fixed $b>0$, take $j=\lfloor bs^2/k\rfloor$. Since
$s\ge2d_R$, we have $j\ge c_bL^{1/4}$ for all sufficiently large $R$,
uniformly in $s$. Consequently,
\begin{equation}
    \label{eq-Pjsbs2}
     \mathrm P_x(J_s>bs^2)\le C_b\exp(-c_bL^{1/4}).
\end{equation}

To cover all radii in the lemma, take $r_j=2^j(2d_R)$ for $j\ge0$
with $r_j\le2R$. If $r_j\le r<2r_j$, then
\[
 B(3r)\subseteq B(6r_j),\qquad \delta r^2\ge\delta r_j^2.
\]
Apply the probability bound \eqref{eq-Pjsbs2} with $s=6r_j$ and $b=\delta/36$,
and take a union bound over the $O(1+L)$ choices of $j$. The polynomial
factor is absorbed by decreasing $c_\delta$, which proves
\eqref{eq-sparse}.
\end{proof}

\section{The walk M(3,2,2)}\label{sec-m322}
Throughout this section, $S$ is the $M(3,2,2)$ walk. We prove
Theorem \ref{thm-m322} by first establishing the bound given in Proposition \ref{prop-cross}. 

\subsection{The expected change before a cube exit}\label{sec-cube}
We begin by expressing the expected change
of a function over one cube exit as a weighted sum of its one-step
changes. The weights will identify the sites where we can ensure a
negative contribution. We then choose the function and show in Lemma \ref{lem-cubedrift} that these negative contributions outweigh those from the remaining sites.

Fix a sufficiently large $R$ throughout this subsection. Write
$\EE_{z_0}$ for expectation under $\PP_{z_0}$ and
$\FF_n=\sigma(S_0,\ldots,S_n)$. Recall that $\tau_s$ is the first time
at which $|S_n|\ge s$.
Let $U_n=\pi(S_n)$, where $\pi(x,y,z)=(x+z,y)$. Given $\FF_n$, the
next increment of $U$ is uniform on the four planar directions. Hence
after every finite $(\FF_n)$-stopping time its future increments,
conditionally on the full history at that time, have the law of planar
simple random walk. In particular $\tau_s<\infty$ almost surely.

Fix a stopping time $T$. On $\{T<\tau_{2R},\ |S_T|\ge2d_R\}$, let
$w=S_T$ and $r=|w|$. We use the cube
\[
 Q=\{u\in\ZZ^3:\|u-w\|_\infty<r/20\},\qquad
 \sigma=\inf\{n\ge T:S_n\notin Q\}.
\]
Let $\mathcal D_T=\{U_n:0\le n<T\}$ be the projected past before
this cube is chosen. For $v\in Q$,
$|\pi v-\pi w|\le\sqrt2\,|v-w|<\sqrt6\,r/20<r$.
Together with $|\pi w|\le\sqrt2\,r$, this gives
\begin{equation}
    \label{eq-piQsetBpir}
     \pi Q\subset B(\pi w,r)\subset B(3r).
\end{equation}
We will use the same disk $B(\pi w,r)$ both to bound visits to one
projected site and to bound the time spent in the projected past.

Write $\PP_T$ and $\EE_T$ for probability and expectation conditional
on $\FF_T$. The expected time before exit is comparable to $r^2$:
\begin{equation}\label{eq-localexit}
  cr^2\le \EE_T(\sigma-T)=\EE_T|S_\sigma-w|^2\le Cr^2.
\end{equation}
The upper bound follows from the definition of $\sigma$. Also note that $(|S_{T+j}-w|^2-j)_{j\geq 0}$ is a martingale since every increment is centered and has squared length one. Optional stopping at $(\sigma-T)\wedge n$ first gives $\EE_T(\sigma-T)<\infty$ and then, by bounded convergence for the squared displacement, the equality
in \eqref{eq-localexit}. At exit at least one coordinate has displacement
of absolute value at least $r/20$, which gives the lower bound.

To compute the expected change up to exit, we group departures by their
starting sites. Let $P_{xy}$ and $P_{yz}$ be the simple random walk
kernels in the $xy$- and $yz$-planes, respectively. For $0\le a\le1$,
write
\[
 P_a=aP_{xy}+(1-a)P_{yz}.
\]
For $u\in Q$, let $g(u)$ and $h(u)$ denote, respectively, the conditional expected number of departures and first departures from $u$ in $[T,\sigma)$, i.e.,
\begin{equation}
    \label{eq-defgh}
    g(u)=\EE_T\sum_{T\le n<\sigma}\ind_{\{S_n=u\}}, \quad h(u)=\EE_T\sum_{T\le n<\sigma}\ind_{\{S_n=u\} \cap \{u \notin \{S_0,S_1,\dots,S_{n-1}\}\}}.
\end{equation}
We denote the fraction of first departures in these expected
counts by
\[
 a(u)=\frac{h(u)}{g(u)}\quad\text{when }g(u)>0.
\]
Write $I$ for the identity operator. For a function
$f$ on the cube and its exit vertices, the conditional expected change
from $u$ is $(P_{xy}-I)f(u)$ on a first departure and $(P_{yz}-I)f(u)$
on a later departure. Summing over the steps before exit gives
\begin{equation}\label{eq-count}
\begin{aligned}
 \EE_T[f(S_\sigma)-f(w)]
   &=\sum_{u\in Q}\bigl[h(u)(P_{xy}-I)f(u)
                       +(g(u)-h(u))(P_{yz}-I)f(u)\bigr]\\
   &=\sum_{\substack{u\in Q\\g(u)>0}}g(u)(P_{a(u)}-I)f(u).
\end{aligned}
\end{equation}
The finiteness needed in \eqref{eq-count} follows from
\eqref{eq-localexit}, since $f$ is bounded on the cube and its exit vertices.
We note that Kalikow \cite{MR628871} used this weighting by expected occupation to
compute exit distributions in a random environment, as described in \cite[Section 4.2.1]{MR2261885}.

The next lemma bounds these fractions. At a site whose projection is
absent from $\mathcal D_T$, the first hit provides a first departure,
whereas the planar Green function bounds the expected total number of
departures after that hit. A possible two-step return also ensures a
uniform upper bound strictly below one at every site.

\begin{lemma}\label{lem-weights}
There is an absolute constant $c_0>0$ such that, for all sufficiently large
$R$ and all $u$ with $g(u)>0$,
\begin{equation}
  a(u)\ge c_0/L\quad\text{if }\pi u\notin\mathcal D_T; \quad a(u)\le a_*:=16/17.
\end{equation}
\end{lemma}
\begin{proof}
Let $\xi_u$ be the first visit to $u$ at or after $T$. On $\{\xi_u<\sigma\}$,
each visit to $u$ during $[\xi_u,\sigma)$ projects to a visit to $\pi u$
before $U$ leaves $B(\pi w,r)$ after $\xi_u$, since this disk contains
$\pi Q$ by \eqref{eq-piQsetBpir}. Conditioning at $\xi_u$ and applying the diagonal bound in
Lemma \ref{lem-green} gives
\begin{equation}
  g(u)\le \PP_T(\xi_u<\sigma)\,G_{B(\pi w,r)}(\pi u,\pi u)
      \le C L\,\PP_T(\xi_u<\sigma).
\end{equation}
If $\pi u\notin\mathcal D_T$, the hit at $\xi_u<\sigma$ gives a first
departure, so $h(u)=\PP_T(\xi_u<\sigma)$. This proves the lower bound.

For the upper bound, a site visited before $T$ has $h(u)=0$. At any other
site, given the first hit in $Q$, choose a $y$ step that stays in the cube.
Such a sign exists since $r/20\ge2$ for all large $R$. This step and its reversal occur with
conditional probability $1/16$, and on this event the walk departs from
$u$ a second time before exit. Hence $g(u)\ge(1+1/16)h(u)$.
\end{proof}

In view of \cite[Example 2.5.9]{MR3587911}, it is natural to choose a positive function for \eqref{eq-count} that is comparable to a negative power of the distance to the origin. However, a radial inverse power alone does not give the required negative change near the shared $y$ axis. Inspired by \cite[the proof of Theorem 1.5]{MR3167134}, we modify it by an angular factor, see Remark \ref{rem-angular-cor} for an explanation of the choice of the factor. With $c_0$ as in Lemma \ref{lem-weights}, set $\varepsilon=c_0/L$. For constants $b,c_1>0$, to be chosen sufficiently small below,
put $\kappa=c_1\sqrt{\varepsilon}$ and define
\begin{equation}\label{eq-potential}
\begin{aligned}
 H(\theta)&=\sqrt{\varepsilon+\theta^2}+\theta^2,\\
 F(v)&=\log\rho+bH(\theta),\\
 V(v)&=e^{-\kappa F(v)},
\end{aligned}
\end{equation}
where $\rho=|v|$ and $\theta=x/\rho$ for $v=(x,y,z)\ne0$. Set $V(0):=1$. The function $V$ depends on $R$ through $\varepsilon$ and $\kappa$. Since $H$ is bounded for $|\theta|\le1$,
\begin{equation}
    \label{eq-Vvkappaasym}
     V(v)\asymp |v|^{-\kappa},\qquad v\ne0,
\end{equation}
with constants independent of $R$. The following bounds provide the
two estimates needed in the weighted sum \eqref{eq-count}.

\begin{lemma}\label{lem-test}
There are fixed $b,c_1>0$ for which, for all large $R$ and all
$|v|\ge3d_R/2$,
\begin{align}
 (P_a-I)V(v)&\le -c\kappa^2 V(v)/|v|^2
       &&(\varepsilon\le a\le a_*),\label{eq-gooddrift}\\
 (P_a-I)V(v)&\le C\kappa^2 V(v)/|v|^2
       &&(0\le a\le a_*).\label{eq-baddrift}
\end{align}
\end{lemma}
\begin{proof}
Under $P_a$, the total probabilities of moving along the three
coordinate axes are $(a,1,1-a)/2$, with equal probabilities for
the two directions along each axis. For a smooth function $f$,
the second-order term in the Taylor expansion of $(P_a-I)f$ is
$(\mathcal L_af)/4$, where
$$
 \mathcal L_a :=a \partial_{xx} +\partial_{yy} +(1-a)\partial_{zz}.
$$
We first reduce the desired bounds on $(P_a-I)V$ to lower bounds for $\mathcal L_aF$.

Since $|\theta|\le1$, one has
$$
 |H'|\le3,\qquad |H''|\le2+\varepsilon^{-1/2},\qquad
 |H'''|\le3/\varepsilon.
$$
Moreover, derivatives of $\theta$ of order $j$ are $O(\rho^{-j})$. Therefore, in a unit neighborhood of $v$ we have $|\nabla F|\le C/\rho$,
and all third derivatives of $F$ are bounded by $C/(\varepsilon\rho^3)$. Let $\xi$ have law $P_a(0,\cdot)$ and put
$\Delta F=F(v+\xi)-F(v)$. Then
$$
 |\Delta F|\le C/\rho,\qquad
 \EE(\Delta F)^2\le C/\rho^2 .
$$
In particular, $\kappa|\Delta F|\le C\kappa/\rho\to0$ as $R \to \infty$.
Using $e^{-t}\le1-t+Ct^2$ for bounded $t$, we obtain
$$
\begin{aligned}
 \frac{(P_a-I)V(v)}{V(v)}
 &=\EE\bigl[e^{-\kappa\Delta F}\bigr]-1\\
 &\le-\kappa\EE[\Delta F]
       +C\kappa^2\EE[(\Delta F)^2]\\
 &\le-\kappa(P_a-I)F(v)+C\kappa^2/\rho^2 .
\end{aligned}
$$
Opposite steps cancel the first-order terms in the Taylor
expansion of $F$, so for large $R$,
$$
 (P_a-I)F(v)
 =\frac14\mathcal L_aF(v)
                  +O((\varepsilon\rho^3)^{-1}).
$$
The error is $o(\sqrt\varepsilon/\rho^2)$ uniformly for
$\rho\ge3d_R/2$ because $\varepsilon^{3/2}d_R\to\infty$ as $R \to \infty$. Since $\kappa=c_1\sqrt{\varepsilon}$, for all sufficiently large $R$,
\begin{equation}\label{eq-local-reduction}
 \frac{(P_a-I)V(v)}{V(v)}
 \le-\frac{\kappa}{4}\mathcal L_aF(v)
   +C\frac{\kappa^2}{\rho^2}.
\end{equation}

We now estimate $\mathcal L_aF$. Write $\omega = (\omega_1, \omega_2, \omega_3) = v/\rho$ and let $C_a=\operatorname{diag}(a,1,1-a)$. Since $|\omega|^2 = 1$, we have $B:=1-\omega^{\mathsf T}C_a\omega=(1-a)\theta^2+a\omega_3^2$.
Differentiating $\log \rho$ and $\theta=x/\rho$ gives
\[
 \rho^2\mathcal L_a\log \rho=2B,\qquad
\rho^2\mathcal L_a\theta=\theta(1-2a-3B),\qquad
\rho\nabla\theta=e_1-\theta\omega.
\]
Put $s=\theta H'(\theta)$ and
$q=(e_1-\theta\omega)^{\mathsf T}C_a(e_1-\theta\omega)$.
Then $0\le s\le|\theta|+2\theta^2\le3$, and $q\ge a(1-\theta^2)^2$ from its first coordinate, and
$H''(\theta)=2+\varepsilon/(\varepsilon+\theta^2)^{3/2}$.
The chain rule gives
\[
 \rho^2\mathcal L_a F
 =(2-3bs)B+bs+b(2q-2as)
       +\frac{b\varepsilon q}{(\varepsilon+\theta^2)^{3/2}} .
\]
Since $(1-\theta^2)^2-s\ge1-|\theta|-4\theta^2\ge-4\theta^2$,
we have $2q-2as\ge-8a\theta^2$. Also
$B\ge(1-a_*)\theta^2$ for $a\leq a_*$. Hence, for $b$ sufficiently small,
$$
 (1-3bs)B+b(2q-2as)
 \ge\bigl[(1-3bs)(1-a_*)-8ba\bigr]\theta^2
 \ge0 .
$$
It follows that
\begin{equation}
 \rho^2\mathcal L_aF
 \ge B+b\left[s+\frac{\varepsilon a(1-\theta^2)^2}
 {(\varepsilon+\theta^2)^{3/2}}\right]\ge0.
\end{equation}
This yields \eqref{eq-baddrift} by \eqref{eq-local-reduction}. Moreover, when $a\ge \varepsilon$, the bracket is at least $c\sqrt\varepsilon$: For $\theta^2\ge\varepsilon$ this follows from $s\ge\theta^2/\sqrt{\varepsilon+\theta^2}\ge\sqrt{\varepsilon/2}$;
for $\theta^2<\varepsilon\le1/4$ the second term is at least
$\varepsilon^2(1-\varepsilon)^2/(2\varepsilon)^{3/2}$. This gives \eqref{eq-gooddrift} for sufficiently small $c_1$.
\end{proof}

\begin{remark}\label{rem-angular-cor}
    On the shared $y$-axis, it is direct to check that the radial term satisfies $\mathcal L_a\log\rho=0$. There the angular correction contributes $baH''(0)/\rho^2$ to $\mathcal L_aF$. Since $H''(0)=2+\varepsilon^{-1/2}$, this contribution is at least $b\sqrt{\varepsilon}/\rho^2$ when $a\ge\varepsilon$. This explains why $H$ in \eqref{eq-potential} is promising and why we let $\kappa$ be a small multiple of $\sqrt{\varepsilon}$ given \eqref{eq-local-reduction}.
\end{remark}

We can now prove the conditional expectation inequality for a cube
exit. In \eqref{eq-count}, sites whose projections are outside
$\mathcal D_T$ contribute negatively by \eqref{eq-gooddrift}.
The contribution from the other sites is controlled by the expected
time spent in $\mathcal D_T$. If this set occupies only a small fraction
of the surrounding disk, Lemma \ref{lem-green} makes that expected time
small compared with the expected exit time in \eqref{eq-localexit}.
We return to the cube centered at $w=S_T$, with $r=|w|$ and exit time
$\sigma$.

\begin{lemma}\label{lem-cubedrift}
There is a fixed $\delta>0$ such that for all sufficiently large $R$,
\begin{equation}\label{eq-cubedrift}
 |\mathcal D_T\cap B(3r)|\le\delta r^2
 \quad\Longrightarrow\quad \EE_T V(S_\sigma)\le V(w).
\end{equation}
\end{lemma}
\begin{proof}
Put $E=\mathcal D_T\cap B(\pi w,r)$. This set is $\FF_T$-measurable
and has at most $\delta r^2$ sites under the condition in
\eqref{eq-cubedrift}, since $B(\pi w,r)\subset B(3r)$.
During $[T,\sigma)$, the projection $\pi(S)$ stays in $B(\pi w,r)$, so its visits
to $\mathcal D_T$ are the same as its visits to $E$. For $0<\delta<1/2$, Lemma \ref{lem-green} therefore gives
\begin{equation}\label{eq-oldtime}
 \EE_T\sum_{T\le n<\sigma}\ind_{\{U_n\in\mathcal D_T\}}
 \le \sum_{v\in E}G_{B(\pi w,r)}(\pi w,v)
 \le C\delta(1+\log(1/\delta))r^2 .
\end{equation}
Here we used that the right-hand side of \eqref{eq-setgreen} is increasing in the set size.

For $u\in Q$ we have $|u|\ge r-\sqrt3\,r/20>3d_R/2$, and
$V(u)/|u|^2\asymp V(w)/r^2$ by \eqref{eq-Vvkappaasym}. In the occupation identity
\eqref{eq-count}, by Lemma \ref{lem-weights}, we can use \eqref{eq-gooddrift} at sites with $\pi u\notin\mathcal D_T$ and \eqref{eq-baddrift} at the remaining sites. By the definition of $g$, see \eqref{eq-defgh},
\[
 \sum_{u\in Q}g(u)=\EE_T(\sigma-T),\qquad
 \sum_{\substack{u\in Q\\\pi u\in\mathcal D_T}}g(u)
   =\EE_T\sum_{T\le n<\sigma}\ind_{\{U_n\in\mathcal D_T\}}.
\]
Subtracting the second sum from the first gives the weight of the sites
with $\pi u\notin\mathcal D_T$. Hence
\[
 \EE_T V(S_\sigma)-V(w)
 \le \frac{\kappa^2 V(w)}{r^2}
 \left[-c\,\EE_T(\sigma-T)
       +C\,\EE_T\sum_{T\le n<\sigma}\ind_{\{U_n\in\mathcal D_T\}}\right].
\]
The lower bound \eqref{eq-localexit} is of order $r^2$, so choosing
$\delta$ small in \eqref{eq-oldtime} proves \eqref{eq-cubedrift}.
\end{proof}

\begin{remark}
(i). A related use of test functions on neighborhoods whose size is proportional to the distance from the origin appears in the recurrence proof for $M(2,1,2)$, see \cite[Lemma 5.2]{Qin2026M212}. \\
(ii). The constants $c_0$, then $b,c_1$, and finally $\delta$ are chosen successively, independently of $R$.
\end{remark}

\subsection{Inward crossings and the proof of Theorem \ref{thm-m322}}\label{sec-finish}
We now iterate Lemma \ref{lem-cubedrift}, continuing only while the
center has distance at least $2d_R$ from the origin and the projected
past satisfies the lemma's sparsity condition. The probability that the sparsity condition fails can be bounded using Lemma \ref{lem-sparse}.

\begin{proposition}\label{prop-cross}
For all sufficiently large $R$ and every starting point $z_0$,
\begin{equation}\label{eq-cross}
 \PP_{z_0}\bigl(\exists\, n\in[\tau_R,\tau_{2R}): |S_n|\le d_R\bigr)
 \le C\exp(-cL^{1/4}),
\end{equation}
where $c,C>0$ are independent of $R$ and $z_0$.
\end{proposition}
\begin{proof}
We may assume that $|z_0|<2R$ since otherwise \eqref{eq-cross} is trivial. Start with $T_0=\tau_R$. At a sampling time $T_j$ ($j\in \NN$),
write $r=|S_{T_j}|$ and stop if any of the following holds:
\[
 T_j\ge\tau_{2R},\qquad r<2d_R,\qquad
 |\mathcal D_{T_j}\cap B(3r)|>\delta r^2.
\]
These conditions correspond, respectively, to leaving the outer ball,
reaching a center close to the origin, and losing the hypothesis of
Lemma \ref{lem-cubedrift}. If none holds, let $T_{j+1}$ be the first
exit from the cube centered at
$S_{T_j}$ with half-side $|S_{T_j}|/20$. We denote the first stopped index by $J$ and set $T_j=T_J$ for $j\geq J$. Each decision is measurable at $T_j$, and $T_J$ is a stopping time, so \eqref{eq-cubedrift} makes $(V(S_{T_j}))_{j\in \NN}$ a nonnegative supermartingale. For each $j<J$, one has $T_j<T_{j+1}<\infty$  by
\eqref{eq-localexit}. Note that $\tau_{2R}<\infty$ as the projection $\pi(S)$ is a planar simple random walk. Thus, both $J$ and $T_J$
are finite almost surely. By Fatou's lemma,
$\EE_{z_0}V(S_{T_J})\le\EE_{z_0}V(S_{\tau_R})\le CR^{-\kappa}$.
On $\{|S_{T_J}|<2d_R\}$ we have $V(S_{T_J})\ge c(2d_R)^{-\kappa}$, so
\begin{equation}\label{eq-smallcenter}
 \PP_{z_0}(|S_{T_J}|<2d_R)
 \le C(2d_R/R)^\kappa\le C e^{-cL^{1/4}}.
\end{equation}
Here we used that $\kappa\asymp L^{-1/2}$.

For each $j<J$, by the stopping rule, $|S_{T_j}|\ge 2d_R$. For $T_j\le n<T_{j+1}$ the walk is still in the
cube centered at $S_{T_j}$, so
\[
 |S_n|\ge |S_{T_j}|-|S_n-S_{T_j}|
       \ge (1-\sqrt{3}/20)|S_{T_j}|>d_R.
\]
Since $T_0=\tau_R$, any $n\in[\tau_R,\tau_{2R})$ with $|S_n|\le d_R$
must therefore satisfy $T_J\le n<\tau_{2R}$. In particular, on the event in \eqref{eq-cross}, one has $T_J<\tau_{2R}$, and if $|S_{T_J}|\ge 2d_R$, then 
\[
 |\mathcal D_{T_J}\cap B(3r)|>\delta r^2,
 \qquad r=|S_{T_J}|\in[2d_R,2R).
\]
As in Lemma \ref{lem-sparse}, let $\mathcal U_R$ be the range of $U=\pi(S)$ before its first exit
from $B(32R)$. For $k<\tau_{2R}$ we have
$|U_k|\le\sqrt{2}\,|S_k|<2\sqrt{2}R$, so
$\mathcal D_{T_J}\subseteq\mathcal U_R$ on $\{T_J<\tau_{2R}\}$.
It follows that
\[
 \PP_{z_0}\bigl(\exists n\in[\tau_R,\tau_{2R}): |S_n|\le d_R\bigr)\le \PP_{z_0}\bigl(|S_{T_J}|<2d_R\bigr)
 +\PP_{z_0}\bigl(\exists r\in[2d_R,2R]:
                   |\mathcal U_R\cap B(3r)|>\delta r^2\bigr).
\]
The two terms on the right-hand side are bounded by \eqref{eq-smallcenter} and
\eqref{eq-sparse}, respectively, proving \eqref{eq-cross}.
\end{proof}

\begin{remark}
The choice of $d_R$ balances two error exponents. More generally,
write $\ell=\log(R/d_R)$. The first-moment estimate in the proof of Lemma \ref{lem-sparse} is $ \EE_x J_s\le C s^2 (1+\ell)L^{-1}$, and the restarting argument (the argument leading to \eqref{eq-Pjsbs2}) gives
\[
 \PP_x(J_s>b s^2)\le C_b\exp\{-c_b L/(1+\ell)\}.
\]
The inward-crossing estimate has exponent $\kappa\ell\asymp \ell/\sqrt L$.
Taking $\ell=L^{3/4}$ makes both exponents of order $L^{1/4}$.
\end{remark}

We are now ready to prove Theorem \ref{thm-m322}.

\begin{proof}[Proof of Theorem \ref{thm-m322}]
 Let $\nu$ be the first visit to zero at or after $\tau_R$. On
$\{\nu<\tau_{2R}\}$, all visits to zero from $\nu$ until $\tau_{2R}$
project to visits of $U=\pi(S)$ to zero before $U$ leaves $B(4R)$ after $\nu$.
Given $\FF_\nu$, the expected number of these visits, including the first,
is at most $C\log R$ by Lemma \ref{lem-green}. Thus \eqref{eq-cross} gives
\[
 \EE_{z_0}\sum_{\tau_R\le n<\tau_{2R}}\ind_{\{S_n=0\}}
 \le C\log R\,\PP_{z_0}(\nu<\tau_{2R})
 \le C\log R\,e^{-c(\log R)^{1/4}} .
\]
For dyadic $R=2^j$, summing these estimates at radii $2^k$, $k\ge j$, yields 
\[
\EE_{z_0} \sum_{n \ge \tau_R }\ind_{\{S_n=0\}} \leq C \sum_{k\geq j} k e^{-c k^{1/4}} \leq C e^{-c'(\log R)^{1/4}},
\]
where we used that $\sum_{k\ge j} k e^{-c k^{1/4}}\le C e^{-c' j^{1/4}}$. Fix $R>|z_0|$. Applying the optional stopping theorem to the martingale $|S_k|^2-k$ at
$\tau_R\wedge n$ gives 
\[
(R+1)^2 \ge \EE_{z_0}|S_{\tau_R\wedge n}|^2=\EE_{z_0}(\tau_R\wedge n)+|z_0|^2 \leq n+|z_0|^2.
\]
Therefore,
\begin{equation}\label{eq-stoptime}
 \PP_{z_0}(\tau_R>n)\le\frac{(R+1)^2}{n},\qquad
 \PP_{z_0}(\tau_R\le n)\le\frac{|z_0|^2+n}{R^2},
\end{equation}
where we also used that $|S_{\tau_R\wedge n}|\ge R$ on $\{\tau_R\le n\}$.
The visits after time $n$ are covered by the visits after $\tau_R$
and, when $n<\tau_R$, the visits between $n$ and $\tau_R$. Thus
\[
 \EE_{z_0}\sum_{k\ge n}\ind_{\{S_k=0\}}
 \le \EE_{z_0}\sum_{k\ge\tau_R}\ind_{\{S_k=0\}}
   +\EE_{z_0}\left[\ind_{\{n<\tau_R\}}
              \sum_{n\le k<\tau_R}\ind_{\{S_k=0\}}\right].
\]
For the second term, condition on $\FF_n$. On $\{n<\tau_R\}$, by \eqref{eq-hitquot} and Lemma \ref{lem-green}, the
expected number of visits before $\tau_R$ is at most
\[
 G_{B(4R)}(U_n,0)\le G_{B(4R)}(0,0)\le C\log R.
\]
Combining this with the first inequality in \eqref{eq-stoptime}, we obtain
\[
 \EE_{z_0}\sum_{k\ge n}\ind_{\{S_k=0\}}
 \le C e^{-c'(\log R)^{1/4}}
       +C(\log R)\frac{(R+1)^2}{n}.
\]
Choose a dyadic radius $R\asymp n^{1/4}$.
This proves \eqref{eq-tail} for large $n$. The remaining times are covered
by increasing $C_v$, with $z_0=-v$.

It remains to establish the escape limit. We first compare $\tau_R$ with
$R^2$. For any fixed $0<\beta<1$, \eqref{eq-stoptime} gives
$\PP_{z_0}(\tau_R\notin [R^{2-\beta},R^{2+\beta}])\le CR^{-\beta}$ for
large $R$. For each fixed $\beta$, Borel--Cantelli therefore gives
$\tau_R\in[R^{2-\beta},R^{2+\beta}]$ eventually along dyadic radii $R=2^j$, $j\geq 1$, almost surely. Applying this along a sequence $\beta\downarrow0$ and
intersecting the resulting full-probability events gives
$\tau_R=R^{2+o(1)}$ almost surely along these radii. By \eqref{eq-cross}, almost surely for all sufficiently large
dyadic $R$ we have $d_R<|S_n|<2R$ throughout $[\tau_R,\tau_{2R})$.
On these same intervals $\log n/\log R\to2$ and $\log|S_n|/\log R\to1$.
This proves \eqref{eq-escape}.
\end{proof}

\section{The walks \texorpdfstring{$M(3,3,2)$ and $M(4,3,2)$}{M(3,3,2) and M(4,3,2)}}\label{sec-transfer}

The proof of Theorem \ref{thm-m3-432} is similar to that of Theorem \ref{thm-m322}. We only highlight the main differences.

\begin{proof}[Proof of Theorem \ref{thm-m3-432}]
For $M(3,3,2)$ put $Z_n=S_n$, and for $M(4,3,2)$ put
\[
 Z_n=\Phi(S_n),\qquad \Phi(x_1,x_2,x_3,x_4)=(x_1,x_3,x_2+x_4).
\]
In both cases, conditional on the original history $\FF_n$, the increment
of $Z$ is uniform on the six neighbors of zero at a first departure of $S$
and on the four neighbors in the $yz$-plane at a later departure. Again,
we translate the target vertex to zero and write $\tau_s=\inf\{n\ge0:|Z_n|\ge s\}$.

Let $\tilde{\pi}(x,y,z):=(y,z)$ and $U_n=\tilde{\pi}(Z_n)$.
Given $\FF_n$, the probability that $U$ moves at the next step is $2/3$
at a first departure of $S$ and $1$ at a later departure. Conditional on
a move, its four directions are equally likely. Thus, starting at any
finite stopping time $T$, the successive moves of $U$ form a planar
simple random walk. Holding steps leave its range unchanged, so
Lemma \ref{lem-sparse} applies to $U$ up to its exit from $B(32R)$.
At each visit of the embedded walk, the conditional mean time to the
next move is at most $3/2$. Hence for a planar disk $A$ containing
$U_T$ and $\rho_A=\inf\{n\ge T:U_n\notin A\}$,
\begin{equation}\label{eq-transfer-green}
 \EE_T\sum_{T\le n<\rho_A}\ind_{\{U_n=v\}}
 \le \frac32 G_A(U_T,v),\qquad v\in A.
\end{equation}

On $\{T<\tau_{2R},\ |Z_T|\ge2d_R\}$, take the cube $Q$ centered at
$w=Z_T$, its exit time $\sigma$, and $r=|w|$ as in Section \ref{sec-cube}.
Write $\mathcal D_T=\{U_n:0\le n<T\}$. Let $g(u)$ be the conditional
expected number of departures of $Z$ from $u$ in $[T,\sigma)$.
Define $h(u)$ in the same way, counting only times when $S$ leaves
its current vertex for the first time. If $\tilde{\pi}(u)\notin\mathcal D_T$,
the first hit $\xi_u$ of $u$ by $Z$ at or after $T$ gives a first
departure of $S$ whenever $\xi_u<\sigma$. In four dimensions, later
visits to $u$ may also give first departures from other vertices of
$\Phi^{-1}(u)$. Since $\tilde{\pi}(Q)\subset B(\tilde{\pi}(w),r)$,
conditioning at $\xi_u$ and using \eqref{eq-transfer-green} and
Lemma \ref{lem-green} similarly as in the proof of Lemma \ref{lem-weights}, we have
\[
 g(u)\le CL\,\PP_T(\xi_u<\sigma)\le CL\,h(u)
 \qquad (\tilde{\pi}(u)\notin\mathcal D_T).
\]
For $g(u)>0$ put $\alpha(u)=h(u)/g(u)$. Let $a_*$ and $c_0$ be as in Lemma \ref{lem-weights}. After decreasing $c_0$ if necessary, we therefore have $2\alpha(u)/3\ge c_0/L$ at sites with $\tilde{\pi}(u)\notin\mathcal D_T$.
Let $P_{xyz}$ be the uniform kernel on the six neighbors in
$\ZZ^3$, which equals $(P_{xy}+P_{xz}+P_{yz})/3$. The occupation identity
\eqref{eq-count} now uses $Q_{\alpha(u)}$, where
$$
 Q_\alpha:=\alpha P_{xyz}+(1-\alpha)P_{yz}=
\frac{\alpha}{3}(P_{xy}+P_{xz})
  +\left(1-\frac{2\alpha}{3}\right)P_{yz}.
$$
Write $J(x,y,z)=(x,z,y)$. Since $J$ interchanges the $xy$- and $xz$-planes and $V\circ J=V$, one has $(P_{xy}V)(Ju)=P_{xz}V(u)$ and $(P_{xz}V)(Ju)=P_{xy}V(u)$ and $(P_{yz}V)(Ju)=P_{yz}V(u)$.
Therefore,
\begin{equation}\label{eq-transfer-kernel}
 (Q_\alpha-I)V(u)=\frac12\big[(P_{2\alpha/3}-I)V(u)
                              +(P_{2\alpha/3}-I)V(Ju)\big].
\end{equation}
Here $P_{2\alpha/3}V(Ju):=(P_{2\alpha/3}V)(Ju)$, which may not equal $P_{2\alpha/3}V(u)$.
We have $2\alpha(u)/3\le2/3<a_*$ wherever $g(u)>0$. Since $|Ju|=|u|$, Lemma \ref{lem-test} (recall $\varepsilon=c_0/L$) and
\eqref{eq-transfer-kernel} give the analogous bounds for
$(Q_{\alpha(u)}-I)V(u)$ with $a=2\alpha(u)/3$:
Inequality \eqref{eq-baddrift} holds wherever $g(u)>0$, and
\eqref{eq-gooddrift} holds if also $\tilde{\pi}(u)\notin\mathcal D_T$. The centered unit steps of $Z$ give the exit-time estimate
\eqref{eq-localexit}. Also, $B(\tilde{\pi}(w),r)\subset B(3r)$, so
\eqref{eq-transfer-green} gives \eqref{eq-oldtime} with an extra
factor of at most $3/2$. Choosing the sparsity constant $\delta$ smaller
if necessary, the proof of Lemma \ref{lem-cubedrift} gives
$\EE_T V(Z_\sigma)\le V(w)$ under the same sparsity condition.
Proposition \ref{prop-cross} now follows for $Z$ by the same iteration.
The occupation and escape estimates \eqref{eq-tail} and \eqref{eq-escape}
follow as in Section \ref{sec-finish}: using \eqref{eq-transfer-green} for visits to the target and the martingale $|Z_n|^2-n$
for \eqref{eq-stoptime}, and applying the dyadic exit-time argument.

In four dimensions, each visit of $S$ to $v$ gives a visit of $Z$ to
$\Phi(v)$, and $|Z_n - \Phi(v)| \le \sqrt{2}\,|S_n - v|$. This gives
\eqref{eq-tail} and the lower bound in \eqref{eq-escape} for $S$. Now let
$\widetilde{\tau}_R = \inf\{n \ge 0 : |S_n - v| \ge R\}$. Fix
$0 < \beta < 1$. Applying the second bound in \eqref{eq-stoptime} to the martingale
$|S_n - v|^2 - n$ and using Borel--Cantelli as in Section \ref{sec-finish}, we
obtain that almost surely $\widetilde{\tau}_{2^j} \ge 2^{j(2-\beta)}$
for all sufficiently large $j$. The lower bound just proved implies
that $|S_n - v| \to \infty$. Hence, for all sufficiently large $n$, the
index $j$ satisfying $2^j \le |S_n - v| < 2^{j+1}$ is large enough that
\[
n \ge \widetilde{\tau}_{2^j} \ge 2^{j(2-\beta)}, \qquad
|S_n - v| < 2 n^{1/(2-\beta)}.
\]
Letting $\beta \downarrow 0$ along a sequence gives the upper bound in \eqref{eq-escape}.
\end{proof}

\section{Two applications of the trace condition}\label{sectrace}
It remains to prove transience for $M(3,1,3)$ and $M(3,2,3)$. We use
the criterion of Peres, Popov and Sousi
\cite[Theorem 1.3]{MR3167134}. Suppose the conditional law
of each increment is chosen by an adapted rule from a finite family of
centered, boundedly supported laws on $\RR^3$. If a fixed matrix $A$
satisfies
\begin{equation}\label{other-trace-condition}
 \operatorname{tr}(A\Sigma A^{\mathsf T})>
       2\lambda_{\max}(A\Sigma A^{\mathsf T})
\end{equation}
for every covariance matrix $\Sigma$ in the family, then the walk is transient.

\begin{proof}[Proof of Proposition \ref{prop-all3models}]
A first-departure law is supported on a line or a plane, so we apply the
criterion after grouping three consecutive steps: Each block has a uniformly positive probability of
containing a later departure, whose covariance is positive definite.

More specifically, fix $d_1\in\{1,2\}$ and sample the walk $S$ of type $M(3,d_1,3)$ at
times $3k$. Conditional on $\FF_{3k}$, let $\eta$ be the expected number
of later departures among the next three steps. The three-step
observation of Camarena, Panizo and Ram\'irez
\cite[argument following Theorem 1]{Camarena2020} gives
\[
 \frac1{2d_1}\le\eta\le3.
\]
To see the lower bound, if either of the first two departures is a later
departure, the block already contains one. Otherwise
$S_{3k}$ and $S_{3k+1}$ are both unvisited before times $3k$ and $3k+1$, respectively; on this $\mathcal{F}_{3k+1}$-measurable event, both steps use the first-departure law, so the second reverses the first with conditional probability $1/(2d_1)$ given $\mathcal{F}_{3k+1}$. A reversal makes the third
departure a later one. Thus the conditional probability that the block
contains a later departure is at least $1/(2d_1)$, which also bounds
$\eta$ from below.

Let $\Pi_{d_1}$ be orthogonal projection onto the first $d_1$ coordinate
axes. Since each step is centered given its past, the cross terms
vanish when computing the conditional covariance of a block. Hence
\[
 \operatorname{Cov}(S_{3k+3}-S_{3k}\mid\FF_{3k})
   =\Sigma(\eta):=(3-\eta)\Pi_{d_1}/d_1+(\eta/3)I.
\]
We check the trace condition first at the two endpoints of this family.
Put $\eta_0=1/(2d_1)$. Both $\Sigma(\eta_0)$ and $\Sigma(3)=I$ are
positive definite, so \cite[Proposition 2.4]{MR3167134} supplies
a single matrix $A$ satisfying \eqref{other-trace-condition} for both.
The trace is linear in the covariance matrix, whereas $\lambda_{\max}$
is convex. Since $\Sigma(\eta)$ is affine in $\eta$, the same strict
inequality holds for every $\eta\in[\eta_0,3]$.

It remains to check that only finitely many
conditional block laws can occur. Every departure in the block occurs
within two steps of $S_{3k}$, so its law is determined by the set
\[
 \{S_n-S_{3k}:0\le n<3k\}
       \cap\{u\in\ZZ^3:|u|\le2\}.
\]
There are finitely many such sets, and each block displacement is
bounded by three. The trace criterion therefore proves transience of
$(S_{3k})$. Since the two intermediate positions are within distance
two of $S_{3k}$, infinitely many visits of $S$ to a fixed vertex would
force infinitely many visits of $(S_{3k})$ to a finite set. Thus $S$ is transient.
\end{proof}

\section*{Acknowledgments}

Shuo Qin is supported by the China Postdoctoral Science Foundation under Grant Numbers 2025M773086 and 2026T190815, and the National Natural Science Foundation of China under Grant Number 12271284. 

\textbf{Use of AI tools.} The function $V$ in \eqref{eq-potential} was designed through interaction with GPT 6. The author also used GPT 6.2 to prepare Figure \ref{fig-cube} and review the paper.

\bibliographystyle{plain}
\bibliography{balanced_excited_refs}
\end{document}